\documentclass[amscd]{amsart}

\newtheorem{Thm}{Theorem}

\newtheorem{Lem}{Lemma}  
\newtheorem{Prop}{Proposition}

\theoremstyle{remark}
\newtheorem{Rem}{Remark}
\newtheorem{Def}{Definition}
\newtheorem{Ex}{Example}

\newcommand{\imbed}{\hookrightarrow}
\newcommand{\iso}{{\widetilde \longrightarrow}}
\newcommand{\To}{\longrightarrow}

\newcommand{\oplusl}{\bigoplus\limits}

\def\square{\hbox{\vrule\vbox{\hrule\phantom{o}\hrule}\vrule}}

\newcommand{\N}{{\mathcal N}}
\renewcommand{\O}{{\mathcal O}}
\newcommand{\F}{{\mathcal F}}
\newcommand{\bO}{{\bf O}}

\newcommand{\Ftil}{\tilde F}

\newcommand{\n}{{\frak n}}
\newcommand{\g}{{\frak g}}
\newcommand{\bu}{\bullet}

\newcommand{\gal}{\check{\ }}

\newcommand{\epf}{\square}

\newcommand{\Jse}{J_\searrow}

\newcommand{\Ine}{I_\nearrow}

\newcommand{\Jne}{J_\nearrow}

\newcommand{\I}{{\mathbb I}}

\newcommand{\Dgood}{{D_{\infty/2}}}

\newcommand{\A}{{\mathcal A}}
\newcommand{\B}{{\mathcal B}}
\newcommand{\Zet}{{\mathbb Z}}
\newcommand{\Ce}{{\mathbb C}}

\begin{document}
\title[Semi-infinite cohomology]{On semi-infinite cohomology of
 finite dimensional algebras}
\author{Roman Bezrukavnikov}
\begin{abstract}
We show that semi-infinite cohomology of a finite dimensional graded
algebra (satisfying some additional requirements) is a particular
case of a general categorical construction. An example of this situation
is provided by small quantum groups at a root of unity.
\end{abstract}
\maketitle

\begin{section}{Introduction}
Semi-infinite cohomology of associative algebras was studied, in particular,
 by S.~Arkhipov
(see \cite{Ar1}, \cite{Ar2}, \cite{Ar3}). Recall that the definition 
of semi-infinite cohomology in \cite{Ar1} works in the following set-up.
We are given an associative graded 
 algebra $A$, two subalgebras $B,N\subset 
A$ such that $A=B\otimes N$ as a vector space, satisfying some additional 
assumptions. In this situation the space of semi-infinite Ext's,
$Ext^{\infty/2 +\bu}(X,Y)$ is defined for $X,Y$ in the bounded derived 
category of graded  $A$-modules. The definition makes use of explicit 
complexes.

 In this note we show that under some additional assumptions
 semi-infinite Ext groups
$Ext^{\infty/2 +\bu}(X,Y)$ has a categorical interpretation. More precisely,
given a category $\A$ and subcategory $\B\subset \A$ one can define
for $X,Y\in \A$ the set 
of morphisms from $X$ to $Y$ "through $\B$"; we denote this space
by $Hom_{\A_\B}(X,Y)$. We then show that if $\A$ is the bounded derived
category of $A$-modules, and $\B$ is the full triangulated subcategory
generated by $B$-projective $A$-modules,  
 then, under certain assumptions
one has 
\begin{equation}\label{eq}
Ext^{\infty/2 +i}(X,Y)=Hom_{\A_\B}(X,Y[i]).
\end{equation}

Notice that the right hand
side of \eqref{eq} makes sense for a wide class of pairs $(A,B)$
(an associative 
 algebra, and a  subalgebra), and $X,Y\in D^b(A-mod)$; in particular we do
not need  $A,B$ to be graded.
Thus one may consider \eqref{eq} as providing
a generalization of the definition  of semi-infinite Ext's to this set up.
However, we should warn the reader that under our working assumptions, but
not in general, $\B$
also equals the full triangulated
 subcategory generated by $B$-injective modules, or by modules
(co)induced from a "complemental" subalgebra $N\subset A$, so one has
at least four different obvious generalizations of the definition of
the right-hand side of \eqref{eq}. 

In fact, a description of semi-infinite cohomology similar to \eqref{eq}
in a general situation (in particular, in the case of enveloping algebras
of infinite-dimensional Lie algebras)
requires additional ideas,
and is the subject of a forthcoming joint work with Arkhipov and Positselskii.

\medskip

An example of the situation considered in this paper is provided by  
a small quantum group at a root of unity \cite{qq},
or by the restricted enveloping algebra of a  simple
Lie algebra in positive characteristic.
Computation of semi-infinite cohomology in the former case is due
to S.~Arkhipov \cite{Ar1} (the answer suggested as a conjecture by 
B.~Feigin). This example was a motivation for the present work.
We informally explain  the relation of our Theorem \ref{teorem}
to the answer for semi-infinite cohomology of small quantum groups
in Remark \ref{remuq} below (we plan to derive it from Theorem \ref{teorem}
elsewhere).
 
\medskip

{\bf Acknowledgements.} 
I thank S.~Arkhipov for stimulating interest, and L.~Positselskii for
helpful comments.  I thank the Clay Mathematical Institute
and NSF grant DMS-0071967 for financial support. 
\end{section}

\begin{section}{Categorical preliminaries:
morphisms through a functor}\label{sec1}
Let $\A$, $\B$ be (small) categories, and $\Phi:\B\to \A$ be a functor.
For $X,Y\in Ob(\A)$ define the set of "morphisms from $X$ to $Y$ through
$\Phi$" 
as  $\pi_0$
of  the category of diagrams 
\begin{equation}\label{dia}
X\To \; \Phi(?) \To \; Y, \ \ \ \ \ \ \ \ \ ?\in \B. 
\end{equation}
 This set
 will be denoted by $Hom_{\A_\Phi}(X,Y)$. Thus elements of $Hom_{\A_\Phi}(X,Y)$
are diagrams of the form \eqref{dia}, with two
diagrams identified if there exists a morphism between them.
Composing the two arrows in \eqref{dia} we get a functorial map
\begin{equation}\label{canmorph}
Hom_{\A_\Phi}(X,Y)\To Hom_\A(X,Y).
\end{equation}

If  $\A$, $\B$ are additive and $\Phi$
 is an additive functor, then addition of
diagrams of the form \eqref{dia} is defined by 
$$(X\overset{f}{\to} \Phi( Z)\overset{g}{\to} Y)+(X
\overset{f'}{\to} \Phi(Z')\overset{g'}{\to} Y)= 
(X\overset{f\times f'}{\To} \Phi(Z\oplus Z') \overset{g\oplus g'}{\To} Y);$$
it induces an abelian group structure on $Hom_{\A_\Phi}(X,Y)$.
 Proposition 3 in \cite{Ma}, VIII.2 shows that  for $Z\in \B$
the tautological map 
$$Hom(X,\Phi (Z)) \otimes _{\Zet} Hom(\Phi(Z),Y)\to Hom_{\A_\Phi}
(X,Y)$$ is compatible with addition.

We have the composition map
$$
Hom_\A(X',X)\times Hom_{\A_\Phi}(X,Y)\times  Hom_\A(Y,Y')
\to Hom_{\A_\Phi}(X',Y);
$$
in particular, for $\A$, $\B$, $\Phi$ additive,
 $Hom_{\A_\Phi}(X,Y)$ is an  $End(X)-End(Y)$ bimodule.

Given $\Phi:\A\to \B$, $\Phi':\A'\to \B'$ and $F:\A\to \A'$, $G:\B\to \B'$
with $F\circ \Phi \cong \Phi'\circ G$ we get for $X,Y\in \A$ a map
\begin{equation}\label{F}
Hom_{\A_\Phi}(X,Y)\to Hom_{\A'_{\B'}}(F(X),F(Y)).
\end{equation}

If  the left adjoint functor  $\Phi^*$ to $\Phi$ is defined on $X$, 
then we have
$$Hom_{\A_\Phi}(X,Y)=Hom_\A (\Phi(\Phi^*(X)),Y),$$
 because in this case the above category
 contracts to the subcategory of diagrams of the form $X\overset{can}{\To}
\Phi(\Phi^*(X)) \to Y$, where $can$ stands for the adjunction morphism.
 If the right adjoint functor $\Phi^!$ is defined on $Y$, then 
$$Hom_{\A_\Phi}(X,Y)=Hom_\A(X, \Phi( \Phi^!(Y)))$$ for similar reasons.
In particular, if $\Phi$ is a full imbedding then 
\eqref{canmorph} is an isomorphism
provided either $X$ or $Y$ lie in the image of $\Phi$.


\medskip
 
In all examples below $\A$ will be a triangulated category,
and  $\Phi:\B\to \A$ will be an imbedding of a (strictly)
 full triangulated subcategory. Given
$\B\subset \A$ we will tacitly assume $\Phi$ to be the
imbedding, and write $Hom_{\A_\B}$
 ("morphisms through $\B$")
instead of  $Hom_{\A_\Phi}$.

\begin{Ex}\label{ex}
 Let $M$ be a Noetherian scheme, and $\A= D^b(Coh_M)$ be the
bounded derived category of coherent sheaves on $M$; let 
$I:\B\imbed \A$
 be the full
  subcategory of complexes whose cohomology ia
supported on a closed subset $i:N\imbed M$. 
Then the functor $I\circ I^!=i_* \circ
i^!$ takes values in a larger derived category of
quasi-coherent sheaves (i.e. ind-coherent sheaves), and $I\circ I^*=
i_* \circ i^*$ takes values
in the Grothendick-Serre dual category,
the derived category of pro-coherent sheaves (introduced
in Deligne's appendix to \cite{H}). Still we have
 $$Hom_{\A_\B}(X,Y)= Hom (X,i_*(i^!(Y)))=Hom (i_*(i^*(X)), Y).$$
In particular, if $X=\O_M$ is the structure sheaf, we get
\begin{equation}\label{cohsup}
Hom_{\A_\B}(\O_M,Y[i])=H_N^i(Y),
\end{equation}
where  $H_N^\bu(Y)$ stands for cohomology
with support on $N$ (local cohomology) \cite{H}.
\end{Ex}

\end{section}

\begin{section}{Recollection of the definition of $Ext^{\infty/2+\bu}$}
\label{assu}
All algebras below will be associative and unital algebras over a field.

\medskip

 We recall a variant of definition of semi-infinite Ext's
 (available under certain restrictions on the 
algebra and subalgebras) suited for our purpose
 (see e.g. \cite{LNM}, \S 2.4, pp 180-183, for this definition
in the particular case of small quantum groups; the general case is 
analogous). 

{\it We make the following assumptions.} 
A $\Zet$-graded algebra $A$  and graded
subalgebras $A^0$, $A^{\leq 0}$, $A^{\geq 0}$
 $\subset A$ are fixed and satisfy the  following 
conditions:

(1)
 $A^{\leq 0}$, $A^{\geq 0}$ are graded by, respectively, 
$\Zet^{\leq 0}$, $\Zet^{\geq 0}$, and $A^0=A^{\leq 0}\cap A^{\geq 0}$ is the
component of degree 0 in $A^{\geq 0}$ and in
$A^{\leq 0}$.

(2) The  maps $A^{\geq 0}\otimes_{A^0} A^{\leq 0} \to A$
and $A^{\leq 0}\otimes_{A^0} A^{\geq 0} \to A$ provided by the 
multiplication map are isomorphisms.

(3)  $A$ is 
 finite dimensional;
$A^0$ is semisimple, and $A^{\geq 0}$ is self-injective (i.e. the free 
$A^{\geq 0}$-module is injective).

\medskip

 By a "module" we will mean a finite dimensional
graded module, unless stated otherwise.
By $A-mod$ we denote the category of (graded finite dimensional)
$A$-modules.

Recall that a bounded below complex of graded modules
is called {\it  convex} if the weights "go down", i.e. for any $n\in \Zet$
the sum of weight spaces of degree more than $n$ is finite dimensional.
A  bounded below complex of graded modules
is called {\it  concave} if the weights "go up" in the similar sense.

\begin{Lem}\label{resolve_b}
i) Any  $A$-module
admits a  right convex resolution by $A$-modules, which are 
injective as $A^{\geq 0}$-modules. 
It also admits a right concave resolution by $A$-modules, which are
$A^{\leq 0}$-injective.

ii) Any finite complex of    $A$-modules
is a quasiisomorphic subcomplex of a bounded
below convex complex of $A^{\geq 0}$-injective $A$-modules. It is also
 a quasiisomorphic subcomplex of a bounded
below concave complex of $A^{\leq 0}$-injective $A$-modules.

\end{Lem}

\proof To deduce (ii) from (i) imbed given finite complex 
$C^\bu\in Com(A-mod)$ into a complex of $A$-injective modules $I^\bu
\in Com^{\geq 0}(A-mod)$ (notice that condition 
(2) above implies that an $A$-injective module is also $A^{\geq 0}$ and
 $A^{\leq 0}$ 
injective), and apply (i) to the module of cocycles 
$Z^n =I^n/d(I^{n-1})$ for large $n$.

To check (i) it suffices to find for any $M\in A-mod$ an imbedding
$M\imbed I$, where $I$ is  $A^{\leq 0}$ injective, and if $n$ is such that
all graded components $M_i$ for $i<n$ vanish, then
$M_n\iso I_n$. (This would prove the second part of the statement; the first
one is obtained from the first one by renotation.)
It suffices to take $I=CoInd_{A^{\geq 0}}^A (Res_{A^{\geq 0}}^A (M))$.
It is indeed $A^{\leq 0}$-injective, because of the equality
\begin{equation}\label{resind}
Res^A_{A^{\leq 0}}(CoInd_{A^{\geq 0}}^A(M))=CoInd_{A^0}^{A^{\leq 0}}(M)), 
\end{equation} 
 which is a consequence of assumption (2) above. \epf

\medskip

We set $D=D^b(A-mod)$.

\begin{Def}\label{defpol}
 (cf. \cite{LNM}, \S 2.4) The assumptions (1--3) are enforced.
Let $X,Y\in D$. Let $\Jse^X$ be
a  convex bounded below complex of 
$A^{\geq 0}$-injective (= projective) modules quasiisomorphic to $X$, and $\Jne ^Y$
be a concave  bounded below complex of 
$A^{\leq 0}$-injective modules quasiisomorphic to $Y$. Then one defines
\begin{equation}\label{defini}
Ext^{\infty/2+i}(X,Y)=H^i(Hom^\bu (\Jse ^X, \Jne^Y)).
\end{equation}

\end{Def}

\begin{Rem}
Independence of  the right-hand side of \eqref{defini}
on the choice of resolutions
$\Jse^X$, $\Jne^Y$ follows from the argument below. Since particular
complexes used in \cite{Ar1} to define $Ext^{\infty/2+\bu}$ satisfy our 
assumptions, we see that this definition agrees with the one in {\it loc.
cit.}
\end{Rem}

\begin{Rem}
Notice that $Hom$ in the right-hand side of \eqref{defini} is $Hom$
in the category of graded modules. As usual, it is often convenient
to denote by $Ext^{\infty/2+i}(X,Y)$ 
the graded space which in present notations
is written down as $\oplusl_n Ext^{\infty/2+i}(X,Y(n))$, where
$(n)$ refers to shift of grading by $-n$.
\end{Rem}

\begin{Rem}
The next standard Lemma shows that conditions on the resolutions
$\Jse^X$, $\Jne^Y$ used in the \eqref{defini} 
 can be formulated in terms of the subalgebra $A^{\geq 0}$ alone
(or, alternatively, in terms of $A^{\leq 0}$ alone);
this conforms with the fact that the left-hand side of \eqref{teoreq}
in Theorem \ref{teorem} below
depends only on $A^{\geq 0}$. However, existence of a "complemental"
subalgebra $A^{\leq 0}$ is used in the construction of a resolution 
$\Jse^X$ with required properties.
\end{Rem}

\begin{Lem}\label{fifi}
An $A$-module is $A^{\leq 0}$-injective iff it is has a filtration with
subquotients of the form $CoInd_{A^{\geq 0}}^A(M)$, $M\in A^{\geq 0}-mod$.
\end{Lem}

\proof The "if" direction follows from semisimplicity of $A^0$,
and  equality \eqref{resind} above.
To show the "only if" part let $M$ be an $A^{\leq 0}$-injective $A$-module.
 Let $M^-$ be its graded
component of minimal degree; then the canonical morphism
\begin{equation}\label{mm}
M\to CoInd_{A^0}^{A^{\leq 0}}(M^-)
\end{equation}
 is surjective. If $M$ is actually
an $A$-module, then the  projection
$M\to M^-$ is a surjection of $A^{\geq 0}$-modules, hence yields
a morphism 
\begin{equation}\label{mmm}
M\to CoInd_{A^{\geq 0}}^A (M^-).
\end{equation}
 \eqref{resind} shows that $Res^A_{A^{\leq 0}}$ sends \eqref{mmm}
into \eqref{mm}; in particular \eqref{mmm} is surjective. Thus the top
quotient of the required filtration is constructed, and the proof is finished
by induction. \epf

\begin{Rem} In two special cases $Ext^{\infty/2+i}(X,Y)$
coincides with a traditional derived functor. First, suppose that
 $Res ^A_{A^{\geq 0}}(X)$ has finite injective (equivalently,
projective) dimension;
 then one can use a finite complex $\Jse^X$ in \eqref{defini} above.
It follows
immediately, that in this case we have $$Ext^{\infty/2+i}(X,Y)\cong
Hom(X,Y[i]).$$

On the other hand, suppose that  $Res ^A_{A^{\leq 0}}(Y)$
has finite injective dimension, so that the complex $\Jne^Y$ in
\eqref{defini} can be chosen 
to be finite. To describe semi-infinite Ext's in this case we need another 
notation. Let $A^*$ denote the co-regular $A$-bimodule; for $M\in A-mod$
let $M\gal=M^*=Hom_A(M,A^*)$ denote the corresponding right $A$-module,
and we use the same notation for the corresponding functor on the derived 
categories. Let also $S:D^b(A-mod)\to D^+(A-mod)$ be given by $S(Y)
=RHom_A(A^*, Y)$. 
Notice that $A^*$ is $A^{\geq 0}$-projective by 
self-injectivity of $A^{\geq 0}$; thus Lemma \ref{fifi} shows that
$Ext^i_A(A^*, N)=0$ for $i>0$ if $N$ is $A^{\leq 0}$-injective. In particular,
$S(Y) \in D^b(A-mod)$ if $Y|_{A^{\leq 0}}$ has finite injective dimension.
We claim that in this case we have 
$$Ext^{\infty/2+i}(X,Y)\cong X\gal \overset{L}{\otimes}_A  S(Y).$$
This isomorphism an immediate consequence of the next Lemma. We also remark
that if $A$ is a Frobenius algebra, then $S\cong Id$.
\end{Rem}

\begin{Lem} Let  $M,N \in A-mod$ be such that $M$ is 
$A^{\geq 0}$-projective, while $N$ is  $A^{\leq 0}$-injective. Then we have

a) $Ext^i_A(M,N)=0$;    $Ext^i_A(A^*,N)=(R^iS)(N)=0$, $Tor_i^A(M\gal ,S(N))=0$
for $i\ne 0$.

b) The natural map
\begin{equation}\label{isomom}
M\gal \otimes_A S(N)= Hom_A(M,A^*)\otimes_A Hom_A(A^*,N)\To Hom_A(M,N)
\end{equation}
is an isomorphism.
\end{Lem}

\proof The first equality in (a) follows from Lemma \ref{fifi}, and the
second one was checked above. Self-injectivity of $A^{\geq 0}$ shows that
$M\gal$ is $A^{\geq 0}$-projective, and a variant of Lemma \ref{fifi}
ensures that it is filtered by modules induced from $A^{\leq 0}$. Thus
it sufficies to show that $S(N)$ is $A^{\leq 0}$-projective. This follows
from isomorphisms
$$Hom_A(A^*, CoInd_{A^{\geq 0}}^A(N_0))= Hom_{A^{\geq 0}} (A^*, N_0)
\cong Hom_{A^{\geq 0}}((A^{\geq 0})^*, N_0)\otimes _{A^0} A^{\leq 0}.$$

Let us now deduce (b) from (a). Notice that (a) implies that both
sides  of \eqref{isomom} 
are exact in $N$ (and also in $M$),
 i.e. send exact sequences $0\to N' \to N\to N'' \to 0$
with $N'$, $N''$ being $A^{\leq 0}$-injective into exact sequences.
Also \eqref{isomom} is evidently an isomorphism for $N=A^*$.
 For any $A^{\leq 0}$-injective $N$ there exists an exact
sequence $$0\to N\to (A^*)^n\overset{\phi}{\To} 
 (A^*)^m$$
 with 
 image and cokernel of $\phi$ being $A^{\leq 0}$-injective.
Thus both sides of \eqref{isomom} turn it into an exact sequence, which
shows that
\eqref{isomom} is an isomorphism for any $A^{\leq 0}$-injective $N$. \epf

\end{section}

\begin{section}{Main result}
\begin{Thm}\label{teorem}
 Let $\Dgood\subset D$ be the full tringulated subcategory
of $D$ generated by $A^{\geq 0}$-injective (=projective) modules. 
For 
 $X,Y\in D^b(A-mod)$ we have a natural
isomorphism
\begin{equation}\label{teoreq}
Hom_{D_{\Dgood}}(X,Y[i])\cong Ext^{\infty/2+i}(X,Y).
\end{equation}
\end{Thm}

\medskip

The proof of  Theorem \ref{teorem} is based on the following

\begin{Lem}\label{resolve_a}
i) Every graded $A^{\geq 0}$-injective $A$-module 
 admits a concave right resolution
consisting of $A$-injective modules.

ii) A finite complex of graded $A^{\geq 0}$-injective $A$-modules is
 quasiisomorphic to 
a concave bounded below complex of $A$-injective modules.

\end{Lem}
\proof (ii) follows from (i) as in the proof of Lemma \ref{resolve_b}.
(Recall that, according to Hilbert,
 if  a bounded below complex of injectives represents an 
object $X\in D^b$ which has finite injective dimension, then for
large $n$ the module of cocycles is injective.)

To prove (i) it is enough for any $A^{\geq 0}$-injective module $M$ to find
an imbedding $M\imbed I$, where $I$ is $A$-injective, and $M_n\iso I_n$
provided $M_i=0$ for $i<n$. (Notice that cokernel of such an imbedding 
is  $A^{\geq 0}$-injective, because $I$ is  $A^{\geq 0}$-injective
by condition (2).) We can take
$I$ to be $CoInd_{A_{\geq 0}}^A (Res _{A_{\geq 0}}^A(M))$. Then
$I$ is indeed injective, because $M$ is $A^{\geq 0}$-injective by 
semi-simplicity of $A^0$, and condition on weights is clearly satisfied.
 \epf

\begin{Prop} a) Let
$\Jse$ be a  convex bounded below complex of 
$A$-modules. 
Let $\Jse^n$ be the $n$-th stupid truncation of $\Jse$
(thus $\Jse^n$ is a quotient complex of $\Jse$).

Let $Z$ be a finite complex of 
$A^{\geq 0}$-injective $A$-modules. Then we have
\begin{equation}\label{lim}
Hom_{D}(X,Z) \iso \varinjlim Hom_{D} (\Jse^n, Z).
\end{equation}
In fact, for $n$ large enough we have 
$$Hom_{D} (X,Z) \iso Hom_{D}  (\Jse^n, Z).$$
\end{Prop}

\proof Let $\Ine$ be a concave bounded below complex of $A$-injective modules
 quasiisomorphic to $Z$ (which exists by Lemma \ref{resolve_a}(ii)).
Then the left-hand side of \eqref{lim} equals $Hom_{Hot}(\Jse, \Ine)$
where $Hot$ stands for the homotopy category of complexes of $A$-modules.  
Conditions on weights of our complexes ensure that
there are only finitely many degrees for which the corresponding 
graded components both in $\Jse$ and $\Ine$ are nonzero; thus
 any morphism
between graded vector spaces $\Jse$, $\Ine$ factors through the  finite
dimensional sum of corresponding graded components. In particular,
  $Hom^\bu (\Jse^n, \Ine) \iso Hom^\bu (\Jse, \Ine)$ for large $n$, and hence
 $$Hom_{D(A-mod)} (\Jse^n, \Ine)=
Hom_{Hot} (\Jse^n, \Ine) \iso Hom_{Hot} (\Jse, \Ine)$$
for large $n$.
\epf

{\it Proof} of the Theorem. We keep notations of Definition
\ref{defpol}.
It follows from the Proposition that
$$      
Hom_{D_\Dgood}(X,Y[i])=\varinjlim_n Hom_D((\Jse^X)^n, Y[i]).
$$     
The right-hand side of \eqref{teoreq} (defined in \eqref{defini})
equals $H^i( Hom^\bu(\Jse^X, \Jne^Y))$. 
Conditions on weights of  $\Jse^X$, $\Jne^Y$
  show that for large $n$ we have
$$   
 Hom^\bu((\Jse^X)^n, \Jne^Y)\iso Hom^\bu(\Jse^X, \Jne^Y) .
$$     
Lemma \ref{fifi} implies that $Ext^i_A(M_1,M_2)=0$ for $i>0$ if $M_1$
is $A^{\geq 0}$-projective, and $M_2$ is $A^{\leq 0}$-injective. Thus
$$ Hom_D((\Jse^X)^n, Y[i])=H^i( Hom^\bu(\Jse^X, \Jne^Y)).$$
The Theorem is proved. \epf

\begin{Rem}\label{remuq}
This remark concerns with the example provided by a  
small quantum group. So let $\g$ be a simple Lie algebra over $\Ce$, 
$q\in \Ce$ be a root of unity of order $l$, and
 let $A=u_q=u_q(\g)$ be the corresponding small quantum group \cite{qq}.
Let $A^{\geq 0}=b_q\subset u_q$ and
$A^{\leq 0}=b_q^-\subset u_q$ be respectively
 the upper and the lower triangular subalgebras.
Then the above conditions (1--3) are satisfied.

Let $\I$ denote the trivial $u_q$-module. The cohomology $Ext^\bu_{u_q}
(\I,\I)$, and the semi-infinite cohomology
$Ext^{\infty/2+\bu}(\I,\I)$ were computed respectively in \cite{GK} and 
\cite{Ar1}. Let us recall the results of these computations.

Assume for simplicity that $l$ is prime to twice 
the maximal multiplicity of an edge in the Dynkin
diagram of $\g$.
Let $\N\subset \g$ be the cone of nilpotent elements, and $\n\subset \N$ be
a maximal nilpotent subalgebra.
 Then the Theorem  of Ginzburg and Kumar 
 asserts that 
\begin{equation}\label{GinK}
Ext^\bu(\I,\I)\cong \O(\N),
\end{equation} 
the algebra of regular functions on $\N$. Also,
 a Theorem of Arkhipov (conjectured by Feigin) asserts that
\begin{equation}\label{Arki}
Ext^{\infty/2+\bu}(\I,\I)\cong H_\n^d(\N,\O),
\end{equation} 
where $d$ is the dimension of $\n$, and $H_\n$ denotes cohomology with
support on $\n$; one also has $H^i_\n(\N,\O)=0$ for $i\ne d$
 (here the choice of $\n$ is assumed to be compatible with the choice
of an upper triangular subalgebra $b_q\subset u_q$ via isomorphism \eqref{GinK}
in a natural sense).

The aim of this remark is to point out a formal similarity 
between \eqref{Arki} and equality \eqref{cohsup} in Example \ref{ex}
above. Namely, the Ginzburg-Kumar isomorphism \eqref{GinK} yields
a functor $F:D^b(u_q-mod)\to Coh(\N)$, $F(X)=Ext^\bu(\I,X)$, such that
$F(\I)=\O_\N$ is the structure sheaf. It is easy to see that 
if $X\in D^b(u_q-mod)$ has finite projective (equivalently, injective)
homological dimension over $b_q$, then the support of $F(X)$
lies in  $\n$
(here by support we mean set-theoretic rather than scheme-theoretic support,
so the coherent sheaf $F(X)$ may be annihilated by some power of the ideal of
 $\n$). Thus if we assume for a moment that the functor $F$ can be lifted 
to a triangulated functor $\Ftil': D^b(u_q-mod)\to D^b(Coh(\N))$,
then \eqref{F} and  Theorem \ref{teorem} would yield a morphism 
from the left-hand side
to the right-hand side of \eqref{Arki}. Here we say that $\Ftil'$ is a
lifting of $F$ if 
$F\cong R\Gamma \circ \Ftil'$, where $R\Gamma(\F)=\oplusl_i H^i(\F)$
for $\F\in D^b(Coh(N))$. 

It is easy to see that
such a functor $\Ftil'$ does not exist. A meaningful version of the argument
is as follows. Let $\bO$ be the differential graded algebra $RHom_{u_q}(\I,\I)$
(thus $\bO$ is a well-defined object of the categroy of differential graded
algebras with inverted quasiisomorphisms); the Ginzburg-Kumar theorem 
\eqref{GinK} shows that the cohomology algebra $H^\bu(\bO)\cong \O(\N)$.
Let  $DGmod(\bO)$ be the 
triangulated category of
differential graded modules over $\bO$  with inverted 
quasiisomorphisms. Let  $D\subset DGmod(\bO)$ be the full subcategory of 
DG-modules
whose cohomology is a finitely generated module over $H^\bu(\bO)=\O(\N)$, and
let $\Dgood\subset D$ be the full 
triangulated subcategory of DG-modules, whose cohomology is a coherent
sheaf on $\N$ supported (set-theoretically) on $\n$.

We have a functor $\Ftil:D^b(u_q-mod)\to D$
 given by $\Ftil:X\mapsto RHom(\I,X)$.
 It is easy to see that $\Ftil$ sends
 complexes of finite homological
dimension over $b_q$ to $\Dgood$; and that
$\Ftil(\I)=\bO$.  Thus, by Theorem \ref{teorem},
\eqref{F} provides a morphism
$$Ext^{\infty/2+\bu}(\I, \I)\To Hom^\bu_{D_\Dgood}(\bO,\bO).$$
One can then show that this morphism is an isomorphism; and also that
the DG-algebra $\bO$ is {\it formal} (quasi-isomorphic to the DG-algebra
$H^\bu(\bO)$ with trivial differential), which implies that 
$$Hom^\bu_{D_\Dgood}(\bO,\bO)\cong H^\bu _\n(\N, \O)$$ (notice that the 
latter isomorphism is not compatible with homological gradings).
This yields the isomorphism \eqref{Arki}.
\end{Rem}

\end{section}


\begin{thebibliography}{111}

\bibitem[Ar1]{Ar1} Arkhipov, S., {\it Semiinfinite cohomology of
quantum groups,} Comm. Math. Phys. {\bf 188} (1997), no. 2, 379--405.

\bibitem[Ar2]{Ar2} Arkhipov, S., {\it Semi-infinite cohomology of
associative algebras and bar duality,} Internat. Math. Res. Notices
{\bf 17} (1997), 833--863. 

\bibitem[Ar3]{Ar3} Arkhipov, S., {\it Semi-infinite cohomology of
quantum groups. II,} in: Topics in quantum groups and finite-type
invariants, 3--42, Amer. Math. Soc. Transl. Ser. 2, 185, Amer. Math. Soc.,
Providence, RI, 1998.


\bibitem[FS]{LNM} Bezrukavnikov, R., Finkelberg, M., Schechtman, V.,
{\it Factorizable sheaves and quantum groups,} Lecture Notes Math. 
1691, Springer Verlag, 1998.

\bibitem[GK]{GK} Ginzburg, V., Kumar, S., {\it
Cohomology of quantum groups at roots of unity,}
Duke Math. J. {\bf 69} (1993), no. 1, 179--198. 

\bibitem[H]{H} Hartshorne, R., {\it Residues and Duality,}
 Lecture Notes Math. 20, Springer Verlag, 1966.

\bibitem[L]{qq} Lusztig, G., {\it Finite-dimensional Hopf algebras 
arising from quantized
universal enveloping algebra,} 
J. Amer. Math. Soc. {\bf 3} (1990), no. 1, 257--296. 

\bibitem[ML]{Ma} Mac Lane, S.,
{\it Categories for the working mathematician,} second edition, 
Springer Verlag, New York, 1998.
\end{thebibliography}
\end{document}